\documentclass[article,12pt]{amsart}

\usepackage{amsfonts}
\usepackage{amsmath}
\usepackage{mathrsfs}
\usepackage{graphicx}
\usepackage{color}
\usepackage{epsfig}
\usepackage{enumerate}

\numberwithin{equation}{section}
\theoremstyle{plain}
\newtheorem{thm}{Theorem}[section]
\newtheorem{rem}{Remark}[section]

\newtheorem{lem}{Lemma}[section]

\newcommand{\dE}{\mathbb{E}}
\newcommand{\dR}{\mathbb{R}}
\newcommand{\dL}{\mathbb{L}}

\newcommand{\cN}{\mathcal{N}}

\newcommand{\rI}{\mathrm{I}}
\newcommand{\cF}{\mathcal{F}}

\newcommand{\cR}{\mathcal{R}}

\newcommand{\veps}{\varepsilon}
\newcommand{\cMc}{\langle M \rangle}

\newcommand{\ind}{\mbox{1}\kern-.25em \mbox{I}}
\font\calcal=cmsy10 scaled\magstep1
\def\build#1_#2^#3{\mathrel{\mathop{\kern 0pt#1}\limits_{#2}^{#3}}}
\def\liml{\build{\longrightarrow}_{}^{{\mbox{\calcal L}}}}
\def\limp{\build{\longrightarrow}_{}^{{\mbox{\calcal P}}}}

\def\videbox{\mathbin{\vbox{\hrule\hbox{\vrule height1.4ex \kern.6em\vrule height1.4ex}\hrule}}}
\def\demend{\hfill $\videbox$\\}

\keywords{Elephant random walk, Martingales, strong law of large numbers, asymptotic normality}

\begin{document}
\title[A martingale approach for the elephant random walk]
{A martingale approach for the elephant random walk \vspace{1ex}}
\author{Bernard Bercu}
\dedicatory{\normalsize University of Bordeaux, France}
\address{Universit\'e de Bordeaux, Institut de Math\'ematiques de Bordeaux,
UMR 5251, 351 Cours de la Lib\'eration, 33405 Talence cedex, France.}
\thanks{}

\begin{abstract}
The purpose of this paper is to establish, via a martingale approach, some refinements on 
the asymptotic behavior of the one-dimensional elephant random walk (ERW). The asymptotic behavior
of the ERW mainly depends on a memory parameter $p$ which lies between zero and one.
This behavior is totally different in the diffusive regime $0 \leq p <3/4$, the critical regime $p=3/4$,
and the superdiffusive regime $3/4<p \leq 1$. Notwithstanding of this trichotomy, we provide
some new results on the almost sure convergence and the asymptotic normality
of the ERW.
\end{abstract}
\maketitle

\ \vspace{-7ex}
\section{Introduction}
\label{S-I}


Random walks with long-memory arose naturally in applied mathematics, theoretical physics, computer sciences
and econometrics. One of them is the so-called elephant random walk (ERW). It is a one-dimensional discrete-time random 
walk on integers, which has a complete memory of its whole history. It was introduced in 2004 by Sch\"utz and 
Trimper \cite{Schutz04} in order to investigate the long-term memory effects in non-Markovian random walks.
It was referred to as the ERW in allusion to the famous saying that elephants can remember where they have been.
\\ \vspace{-1ex}\par
A wide range of literature is available on the asymptotic behavior of the ERW and its extensions 
\cite{Baur16},\cite{Boyer14},\cite{Col17},\cite{Cressoni13},\cite{Cressoni07},\cite{Da13},\cite{Kumar10},\cite{Kursten16}. 
However, many things remain to be done. The goal of this paper is to answer to several natural questions on the ERW. 
\\ \vspace{-1ex}\par
One of them concerns the influence of the memory parameter $p$ on the almost sure asymptotic behavior of the ERW. 
Quite recently, Baur and Bertoin \cite{Baur16} and independently Coletti, Gava and Sch\"utz \cite{Col17} 
have studied the asymptotic normality of the ERW in the diffusive regime $p<3/4$ as well as in the critical regime $p=3/4$. 
However, very few results are available on the almost sure asymptotic behavior of the ERW in the regime $p\leq 3/4$. 
We shall fill the gap by proving a quadratic strong law as well as a law of iterated logarithm for
the ERW in the regime $p\leq 3/4$.
\\ \vspace{-1ex}\par
Another key question concerns the limiting distribution of the ERW in the superdiffusive regime $p> 3/4$.
Initially, it was suggested by Sch\"utz and Trimper \cite{Schutz04} that, even in the superdiffusive regime,
the ERW has a Gaussian limiting distribution. Later, it was conjectured by Da Silva et al. \cite{Da13}
that this limiting distribution is not Gaussian, see also \cite{Col17},\cite{Par06}.
One can observe that the analytical study of \cite{Da13} is not sufficient to prove the non-gaussianity of the
limiting distribution. We shall provide a rigorous mathematical proof that the limiting distribution is not Gaussian.
Starting from the symmetric initial condition, we will also show that this limiting distribution is sub-Gaussian.
\\ \vspace{-1ex}\par
Baur and Bertoin \cite{Baur16} extensively used the connection to P\'{o}lya-type urns \cite{Harris15}
as well as two functional limit theorems for multitype branching processes due to Janson \cite{Janson04}, see also \cite{Chauvin11}.
Our strategy is totally different as it relies on the theory of martingales. To be more precise, we shall make use of the strong
law of large numbers and the central limit theorem for martingales \cite{Duflo97},\cite{HallHeyde80} as well as 
the law of iterated logarithm for martingales \cite{Stout70},\cite{Stout74}. We strongly believe that our approach could be successfully
extended to ERW with stops \cite{Cressoni13},\cite{Harbola14}, to amnesiac ERW \cite{Cressoni07}, as well as to multi-dimensional ERW
\cite{Cressoni13},\cite{Lyu17}.
\\ \vspace{-1ex}\par
The paper is organized as follows. In Section \ref{S-E}, we introduce the exact ERW and the martingale
we shall extensively make use of. Section \ref{S-MR} is devoted to the main results of the paper.
We establish the almost sure asymptotic behavior as well as the asymptotic normality of the ERW
in the diffusive and critical regimes. Moreover, in the superdiffusive regime, 
we provide the first rigorous mathematical proof that the limiting distribution of the ERW is not Gaussian.
Our martingale approach is described in Appendix A, while all technical proofs of Section \ref{S-MR}
are postponed to Appendices B and C.

\ \vspace{-2ex} \\
\section{The elephant random walk}
\label{S-E}


The one-dimensional ERW is defined as follows. The random walk starts at the origin at time zero,
$S_0=0$. At time $n=1$, the elephant moves to the right with probability $q$ and to the left with probability
$1-q$ where $q$ lies between zero and one. Hence, the position of the elephant at time $n=1$ is given by
$S_1=X_1$ where $X_1$ has a Rademacher $\cR(q)$ distribution.
Afterwards, at any time $n \geq 1$, we choose uniformly at random an interger $k$ among the previous times
$1,\ldots,n$, and we define
\begin{equation*}
   X_{n+1} = \left \{ \begin{array}{ccc}
    +X_{k} &\text{ with probability } & p, \vspace{2ex}\\
    -X_{k} &\text{ with probability } & 1-p ,
   \end{array} \nonumber \right.
   \vspace{2ex}
\end{equation*}
where the parameter $p \in [0,1]$ is the memory of the ERW. 
Then, the position of the ERW is given by
\begin{equation}
\label{POSERW}
S_{n+1}=S_{n}+X_{n+1}.
\end{equation}
\\ \vspace{-3ex}\par
In order to understand well how the elephant moves, it is straightforward to see that 
for any time $n \geq 1$, $X_{n+1} = \alpha_{n} X_{\beta_n}$
where $\alpha_n$ and $\beta_n$ are two independent discrete random variables where  
$\alpha_n$ has a Rademacher $\cR(p)$ distribution while $\beta_n$ is uniformly distributed over
the integers $\{1,\cdots,n\}$. Moreover, $\alpha_n$ is independent of $X_1,\ldots,X_n$.
\\ \vspace{-1ex}\par
Let $(\cF_n)$ be the increasing sequence of $\sigma$-algebras,
$\cF_n=\sigma(X_1,\ldots,X_n)$.
For any time $n \geq 1$, we clearly have 
\begin{equation}
\label{CEX}
   \dE[X_{n+1} | \cF_{n}] = \dE[\alpha_{n}] \times \dE[X_{\beta_n} | \cF_{n}]=(2p-1)\frac{S_n}{n}  \hspace{1cm} \text{a.s.}
\end{equation}
which, together with \eqref{POSERW}, implies that
\begin{equation}
\label{CES1}
   \dE[S_{n+1} | \cF_{n}] = \gamma_{n} S_{n} \hspace{1cm} \text{where} \hspace{1cm}
   \gamma_n=\Bigl( \frac{n+2p-1}{n} \Bigr).
\end{equation}
Moreover,
$$
\prod_{k=1}^n \gamma_k =\frac{\Gamma(n+2p)}{\Gamma(n+1)\Gamma(2p)}
$$
where $\Gamma$ stands for the Euler gamma function. Therefore, let $(M_n)$ be the sequence of random variables defined, 
for all $n \geq 0$, by $M_n = a_n S_n$ where $a_1=1$ and, for all $n\geq 2$,
\begin{equation}
\label{DEFAN}
a_n= \prod_{k=1}^{n-1} \gamma_k^{-1} = \frac{\Gamma(n)\Gamma(2p)}{ \Gamma(n+2p-1)}.
\end{equation}
Since $a_n = \gamma_n a_{n+1}$, we clearly deduce from \eqref{CES1} that for any time $n \geq 1$,
\begin{equation*}
   \dE[ M_{n+1} | \cF_{n}] = M_{n} \hspace{1cm} \text{a.s.} 
\end{equation*}
In other words, the sequence $(M_n)$ is a multiplicative real martingale. Our strategy is to make use of
the martingale $(M_n)$ in order to deduce the asymptotic behavior of $(S_n)$.


\section{Main results.}
\setcounter{equation}{0}
\label{S-MR}


\subsection{The diffusive regime}

Our first result concerns the almost sure convergence of the ERW in the diffusive regime where $0 \leq p <3/4$.

\begin{thm}
\label{T-ASCVG-DR}
We have the almost sure convergence
\begin{equation}
\label{T-ASCVG-DR1}
 \lim_{n \rightarrow \infty} \frac{S_n}{n}=0 \hspace{1cm} \text{a.s.}
\end{equation}
\end{thm}

\noindent
We focus our attention on the almost sure rates of convergence of the ERW.
\begin{thm}
\label{T-ASCVGRATES-DR}
We have the quadratic strong law
\begin{equation}
\label{T-ASCVG-DR2}
 \limsup_{n \rightarrow \infty} \frac{1}{\log n} \sum_{k=1}^n  \Bigl(\frac{S_k}{k}\Bigr)^2=\frac{1}{3-4p} \hspace{1cm} \text{a.s.}
\end{equation}
In addition, we also have the law of iterated logarithm
\begin{eqnarray}
 \limsup_{n \rightarrow \infty} \Bigl(\frac{1}{2 n \log \log n}\Bigr)^{1/2} S_n & = & 
 -\liminf_{n \rightarrow \infty} \Bigl(\frac{1}{2 n \log \log n}\Bigr)^{1/2} S_n \nonumber \\
 & = & \frac{1}{\sqrt{3-4p}} \hspace{1cm} \text{a.s.}
 \label{T-ASCVG-DR3}
\end{eqnarray}
In particular,
\begin{equation}
 \limsup_{n \rightarrow \infty} \frac{S_n^2}{2 n \log \log n}=  \frac{1}{3-4p} \hspace{1cm} \text{a.s.}
 \label{T-ASCVG-DR4}
\end{equation}
\end{thm}

\noindent
Our next result is devoted to the asymptotic normality of the ERW in the diffusive regime $0 \leq p <3/4$.

\begin{thm}
\label{T-AN-DR}
We have the asymptotic normality
\begin{equation}
\label{T-AN-DR1}
\frac{S_n}{\sqrt{n}} \liml \cN \Bigl(0, \frac{1}{3-4p} \Bigr).
\end{equation}
\end{thm}

\begin{rem} One can observe that the additional term $(2q-1)n^{2p-1}/\Gamma(2p)$
is useless in Theorem 2 of \cite{Col17}. Moreover, in the particular case $p=1/2$, one find again
\begin{equation*}
\frac{S_n}{\sqrt{n}} \liml \cN (0,1).
\end{equation*}
\end{rem}

\subsection{The critical regime}
Hereafter, we investigate the critical regime where the memory parameter $p=3/4$.

\begin{thm}
\label{T-ASCVG-CR}
We have the almost sure convergence
\begin{equation}
\label{T-ASCVG-CR1}
 \lim_{n \rightarrow \infty} \frac{S_n}{\sqrt{n} \log n}=0 \hspace{1cm} \text{a.s.}
\end{equation}
\end{thm}

\noindent
The almost sure rates of convergence of the ERW are as follows
\begin{thm}
\label{T-ASCVGRATES-CR}
We have the quadratic strong law
\begin{equation}
\label{T-ASCVG-CR2}
 \limsup_{n \rightarrow \infty} \frac{1}{\log \log n} \sum_{k=2}^n  \Bigl(\frac{S_k}{k \log k}\Bigr)^2=1 \hspace{1cm} \text{a.s.}
\end{equation}
In addition, we also have the law of iterated logarithm
\begin{eqnarray}
 \limsup_{n \rightarrow \infty} \Bigl(\frac{1}{2 n \log n \log \log \log n}\Bigr)^{1/2} S_n & = & 
 -\liminf_{n \rightarrow \infty} \Bigl(\frac{1}{2 n \log n \log \log \log n}\Bigr)^{1/2} S_n \nonumber \\
 & = & 1 \hspace{1cm} \text{a.s.}
 \label{T-ASCVG-CR3}
\end{eqnarray}
In particular,
\begin{equation}
 \limsup_{n \rightarrow \infty} \frac{S_n^2}{2 n \log n \log \log \log n}=  1 \hspace{1cm} \text{a.s.}
 \label{T-ASCVG-CR4}
\end{equation}
\end{thm}

\noindent
One can observe a very unusual rate of convergence in the law of iterated logarithm.
Our next result deals with the asymptotic normality of the ERW in the critical regime $p=3/4$.

\begin{thm}
\label{T-AN-CR}
We have the asymptotic normality
\begin{equation}
\label{T-AN-CR1}
\frac{S_n}{\sqrt{n \log n}} \liml \cN (0,1).
\end{equation}
\end{thm}

\begin{rem} As before, the additional term $(2q-1)\sqrt{n}/\Gamma(3/2)$
is useless in Theorem 2 of \cite{Col17}. 
\end{rem}

\subsection{The superdiffusive regime}
Finally, we focus our attention on the more complicated superdiffusive regime where 
$3/4 < p \leq 1$.

\begin{thm}
\label{T-ASCVG-SR}
We have the almost sure convergence
\begin{equation}
\label{T-ASCVG-SR1}
 \lim_{n \rightarrow \infty} \frac{S_n}{n^{2p-1}}=L \hspace{1cm} \text{a.s.}
\end{equation}
where $L$ is a non-degenrate random variable. This convergence 
also holds in $\dL^4$, which means that
\begin{equation}
\label{T-ASCVG-SR2}
 \lim_{n \rightarrow \infty} \dE\Bigl[ \Bigl| \frac{S_n}{n^{2p-1}} -L \Bigr|^4 \Bigr]=0.
\end{equation}
\end{thm}

\ \vspace{-2ex} \\
\begin{rem}
One can observe that the first three moments of $S_n$ where previously calculated in \cite{Da13} in the
special case $q=1$. However, the analytical study of \cite{Da13} is not sufficient to evaluate the
moments of $L$. 
\end{rem}

\begin{thm}
\label{T-MOM-SR}
The first four moments of $L$ are given by 
\begin{eqnarray}
\label{MOML1}
\dE[L]  &=& \frac{2q-1}{\Gamma(2p)},\\
\label{MOML2}
\dE[L^2] &=& \frac{1}{(4p-3)\Gamma(2(2p-1))}, \\
\label{MOML3}
\dE[L^3] &=& \frac{2p(2q-1)}{(2p-1)(4p-3)\Gamma(3(2p-1))}, \\
\label{MOML4}
\dE[L^4] &=& \frac{6(8p^2-4p-1)}{(8p-5)(4p-3)^2\Gamma(4(2p-1))}.
\end{eqnarray}
\end{thm}

\ \vspace{-2ex} \\
\begin{rem} Our last result provides the first rigorous mathematical proof that, in the
superdiffusive regime, the limiting distribution $L$ of the ERW is not Gaussian.
As a matter of fact, denote by $\mu$ and $\sigma^2$ the mean value and the 
variance of $L$, $\mu=\dE[L]=2q-1$ and $\sigma^2=\dE[(L-\mu)^2]$.
Moreover, let $\alpha$ and $\kappa$ be the skewness and the kurtosis of $L$, 
respectively defined by
$$
\alpha=\frac{\dE[(L-\mu)^3]}{\sigma^3}
\hspace{1cm}\text{and}\hspace{1cm}
\kappa=\frac{\dE[(L-\mu)^4]}{\sigma^4}.
$$
In the special case $q=1/2$, the skewness $\alpha=0$, while the kurtosis
$$
\kappa=\frac{6(8p^2-4p-1)(\Gamma(2(2p-1)))^2}{(8p-5)\Gamma(4(2p-1))}.
$$
It is not hard to see that $\kappa$ is a decreasing function of $p$ such that, 
for all $3/4 < p \leq 1$, $1 \leq \kappa < 3$. It means that $L$ has a sub-Gaussian distribution.
Furthermore, if $q=1$, $\alpha>0$ and if $q=0$, $\alpha<0$.
Finally, the kurtosis $\kappa$ shares the same value for $q=1$ or $q=0$, which can be smaller, greater or
equal to $3$.
\end{rem}

\section*{Appendix A \\ The martingale approach}
\renewcommand{\thesection}{\Alph{section}}
\renewcommand{\theequation}{\thesection.\arabic{equation}}
\setcounter{section}{1}
\setcounter{equation}{0}

We already saw that the sequence $(M_n)$ given, for all $n \geq 0$, by $M_n=a_nS_n$, is a multiplcative real martingale.
Moreover, for any $n \geq 1$, $X_n$ is a binary random variables taking values in $\{+1,-1\}$.
Consequently, $|S_n| \leq n$, which implies that $(M_n)$ is locally square integrable.
The martingale $(M_n)$ can be rewritten in the additive form
\begin{equation}
\label{DECMN}
M_n=\sum_{k=1}^n a_{k} \veps_{k}
\end{equation}
since its increments $\Delta M_n= M_n-M_{n-1}$ satisfy
$\Delta M_{n}= a_{n}S_{n}-a_{n-1}S_{n-1}=a_{n} \veps_{n}$
where $\veps_{n}=S_{n}-\gamma_{n-1}S_{n-1}$. The predictable quadratic variation \cite{Duflo97}
associated with 
$(M_n)$ is given by $\cMc_0=0$ and, for all $n\geq 1$, 
\begin{equation}
\label{DEFIP}
\cMc_n=\sum_{k=1}^n \dE[ \Delta M_k^2|\cF_{k-1}].
\end{equation}
We immediately obtain from \eqref{CES1} that
$\dE[\veps_{n+1} | \cF_{n}]=0$. Moreover, it follows from \eqref{POSERW} together with \eqref{CEX} that
\begin{equation}
\label{CES2}
   \dE[S_{n+1}^2 | \cF_{n}] = \dE[ S_n^2 +2 S_n X_{n+1} +1 | \cF_{n}]=1+ ( 2 \gamma_n -1)S_{n}^2 
   \hspace{1cm} \text{a.s.}
\end{equation}
Consequently, as $\dE[\veps_{n+1}^2 | \cF_{n}]= \dE[ S_{n+1}^2 | \cF_{n}] - \gamma_n^2S_n^2$,
we deduce from \eqref{CES2} that, for all $n\geq 1$, 
\begin{eqnarray}
   \dE[\veps_{n+1}^2 | \cF_{n}]  
   & = &  1 + ( 2 \gamma_n -1)S_{n}^2  - \gamma_n^2S_n^2 =  1 - (\gamma_n -1)^2 S_n^2  \hspace{1cm} \text{a.s.} \nonumber\\
   & = & 1 - (2p-1)^2 \Bigl(\frac{S_n}{n}\Bigr)^2  \hspace{1cm} \text{a.s.}
   \label{CEEPS2}
\end{eqnarray}
By the same token, 
\begin{eqnarray}
   \dE[\veps_{n+1}^4 | \cF_{n}]  
   & = &  1 -3(\gamma_n -1)^4S_n^4  +2 (\gamma_n -1)^2 S_n^2  \hspace{1cm} \text{a.s.} \nonumber\\
   & = & 1 - 3(2p-1)^4 \Bigl(\frac{S_n}{n}\Bigr)^4 + 2(2p-1)^2 \Bigl(\frac{S_n}{n}\Bigr)^2  \hspace{1cm} \text{a.s.}
   \label{CEEPS4}
\end{eqnarray}
On the one hand, if $p=1/2$, $\dE[\veps_{n+1}^2 | \cF_{n}]=1$ and $\dE[\veps_{n+1}^4 | \cF_{n}]=1$ a.s.
On the other hand, we obtain from \eqref{CEEPS2} and \eqref{CEEPS4}
the almost sure upper bounds
\begin{equation}
\label{UBCEEPS}
\sup_{n \geq 0} \dE[\veps_{n+1}^2 | \cF_{n}] \leq 1 \hspace{1cm} \text{and} \hspace{1cm}
\sup_{n \geq 0} \dE[\veps_{n+1}^4 | \cF_{n}] \leq \frac{4}{3}.
\end{equation}
Hereafter, we deduce from \eqref{DECMN}, \eqref{DEFIP} and \eqref{CEEPS2} that 
\begin{equation}
\label{CALIP}
\cMc_n=\sum_{k=1}^n a_k^2 - (2p-1)^2 \zeta_n
\hspace{1cm} \text{where} 
\hspace{1cm} 
\zeta_n = \sum_{k=1}^{n-1} a_{k+1}^2 \Bigl(\frac{S_k}{k}\Bigr)^2.
\end{equation}
The asymptotic behavior of the martingale $(M_n)$ is closely related to the one of 
$$
v_n=\sum_{k=1}^n a_k^2=  \sum_{k=1}^{n}  \Bigl( \frac{\Gamma(k) \Gamma(2p) }{\Gamma(k+2p-1)} \Bigr)^2.
$$
Via standard results on the asymptotic behavior of the Euler gamma function, we have three regimes.
In the diffusive regime where $0 \leq p <3/4$,
\begin{equation}
\label{CVGVNDR}
\lim_{n\rightarrow \infty} \frac{v_{n}}{n^{3-4p}} = \ell \hspace{1cm} \text{where} \hspace{1cm}
\ell=\frac{(\Gamma(2p))^2}{3-4p}.
\end{equation}
In the critical regime where $p=3/4$,
\begin{equation}
\label{CVGVNCR}
\lim_{n\rightarrow \infty} \frac{v_{n}}{\log n} = \frac{\pi}{4}.
\end{equation}
In the superdiffusive regime where 
$3/4 < p \leq 1$, $v_n$ converges to the finite value 
\begin{equation}
\label{CVGVNSR}
\lim_{n\rightarrow \infty} v_n = \sum_{k=0}^{\infty}  \Bigl( \frac{\Gamma(k+1) \Gamma(2p) }{\Gamma(k+2p)} \Bigr)^2
= \sum_{k=0}^{\infty} \frac{(1)_k\,(1)_k\,(1)_k} {(2p)_k\, (2p)_k\,k!} = {}_{3}F_2\Bigl( \begin{matrix}
{1,1,1}\\
{2p,2p}\end{matrix} \Bigl| 1\Bigr)
\end{equation}
where, for any $a\in \dR$, $(a)_k=a(a+1)\cdots(a+k-1)$ for $k\geq 1$, $(a)_0=1$ stands for the Pochhammer symbol and ${}_{3}F_2$ is the
generalized hypergeometric function defined by
\begin{eqnarray*}
{}_{3}F_2 \Bigl( \begin{matrix}
{a,b,c}\\
{d,e}\end{matrix} \Bigl|
{\displaystyle z}\Bigr)
=\sum_{k=0}^{\infty}
\frac{(a)_k\,(b)_k\,(c)_k}
{(d)_k\,(e)_k\, k!} z^k.
\end{eqnarray*}

\section*{Appendix B \\ Proofs of the almost sure convergence results}
\renewcommand{\thesection}{\Alph{section}}
\renewcommand{\theequation}{\thesection.\arabic{equation}}
\setcounter{section}{2}
\setcounter{equation}{0}

\vspace{-3ex}
\subsection*{}
\begin{center}
{\bf B.1. The diffusive regime.}
\end{center}
\ \vspace{-4ex}\\

\noindent{\bf Proof of Theorem \ref{T-ASCVG-DR}.}
First of all, we focus our attention on the proof of the almost sure convergence
\eqref{T-ASCVG-DR1}. We already saw from \eqref{CALIP} that $\cMc_n \leq v_n$.
Moreover, it follows from \eqref{CVGVNDR} that, in the diffusive regime, $v_n$ 
increases to infinity with an arithmetic speed $n^{3-4p}$.
Then, we obtain from the strong law of large numbers for martingales given e.g. by Theorem 1.3.24 of \cite{Duflo97} that
\begin{equation*}
\vspace{-1ex}
\lim_{n\rightarrow \infty} \frac{M_{n}}{v_{n}} = 0 \hspace{1cm} \text{a.s.}
\end{equation*}
which implies that
\begin{equation}
\label{CVGMN}
\lim_{n\rightarrow \infty} \frac{M_{n}}{n^{3-4p}} = 0 \hspace{1cm} \text{a.s.}
\end{equation}
This convergence is not sharp enough to prove \eqref{T-ASCVG-DR1}.
However, thanks to the last part of Theorem 1.3.24,
we also have the almost sure rate of convergence in \eqref{CVGMN}
\begin{equation*}
\label{CVGRATEMN}
\vspace{-1ex}
\frac{M_{n}^2}{v_n}  = O(\log v_n) \hspace{1cm} \text{a.s.}
\end{equation*}
which ensures that
\begin{equation}
\label{CVGRATEMN}
\frac{M_{n}^2}{n^{3-4p}}  = O(\log n) \hspace{1cm} \text{a.s.}
\end{equation}
Hereafter, as $M_n=a_nS_n$, it clearly follows from \eqref{CVGRATEMN} that
\begin{equation*}
\frac{S_{n}^2}{n}  = O(\log n) \hspace{1cm} \text{a.s.}
\end{equation*}
which immediately leads to \eqref{T-ASCVG-DR1}.
\demend


\noindent{\bf Proof of Theorem \ref{T-ASCVGRATES-DR}.}
Denote by $f_{n}$ the explosion coefficient associated with the martingale 
$(M_{n})$ given, for all $n \geq 1$, by
\begin{equation*}
f_{n} = \frac{a_n^2}{v_n}.
\end{equation*}
We clearly obtain from \eqref{CVGVNDR} that $f_n$ converges to zero. Moreover, it follows from
the almost sure convergence \eqref{T-ASCVG-DR1} together with \eqref{CEEPS2} that
\begin{equation*}
\lim_{n \rightarrow \infty} \dE[\veps_{n+1}^2 | \cF_{n}]=1 \hspace{1cm} \text{a.s.}
\end{equation*}
Furthermore, we already saw in \eqref{UBCEEPS} that
\begin{equation*}
\sup_{n \geq 0} \dE[\veps_{n+1}^4 | \cF_{n}] \leq \frac{4}{3} \hspace{1cm} \text{a.s.}
\end{equation*}
Consequently, by virtue of the quadratic strong law for martingales given e.g. in Theorem 3 of \cite{Bercu04},
\begin{equation*}
\lim_{n\rightarrow \infty} 
\frac{1}{\log v_n} \sum_{k=1}^{n} f_{k} \Bigl( \frac{M_{k}^2}{v_k} \Bigr) = 1 \hspace{1cm} \text{a.s.}
\end{equation*}
which implies, via \eqref{CVGVNDR}, that
\begin{equation}
\label{LFQ-DR1}
\lim_{n\rightarrow \infty} 
\frac{1}{\log n} \sum_{k=1}^{n}   \frac{a_k^2 M_{k}^2}{v_k^2} = (3-4p) \hspace{1cm} \text{a.s.}
\end{equation}
Therefore, as $M_n=a_nS_n$ and $n^2 a_n^4$ is equivalent to $(3-4p)^2v_n^2$, we find
from \eqref{LFQ-DR1} that
\begin{equation}
\label{LFQ-DR2}
\lim_{n\rightarrow \infty} 
\frac{1}{\log n} \sum_{k=1}^n  \Bigl(\frac{S_k}{k}\Bigr)^2=\frac{1}{3-4p} \hspace{1cm} \text{a.s.}
\end{equation}
which completes the proof the quadratic strong law \eqref{T-ASCVG-DR2}. We shall now proceed to the proof of the law of iterated logarithm
given by \eqref{T-ASCVG-DR3}. In order to apply the law of iterated logarithm for martingales due to
Stout \cite{Stout74}, see also Corollary 6.4.25 in \cite{Duflo97}, it is only necessary to verify that
\begin{equation}
\label{CONDLIL}
\sum_{n=1}^{+\infty} \frac{a_{n}^4}{v_n^2} < +\infty. 
\end{equation}
This is clearly satisfied since $a_n^4v_n^{-2}$ is equivalent to $(3-4p)^2n^{-2}$ 
and, as is well-known,
\begin{equation*}
\sum_{n=1}^{+\infty} \frac{1}{n^2} =\frac{\pi^2}{6}.
\end{equation*}
Hence, we find from the law of iterated logarithm for martingales that
\begin{eqnarray}
 \limsup_{n \rightarrow \infty} \Bigl(\frac{1}{2 v_n \log \log v_n}\Bigr)^{1/2} M_n & = & 
 -\liminf_{n \rightarrow \infty} \Bigl(\frac{1}{2 v_n \log \log v_n}\Bigr)^{1/2} M_n \nonumber \\
 & = & 1 \hspace{1cm} \text{a.s.}
 \label{LIL-MG-DR}
\end{eqnarray}
As previously seen, the identity $M_n\!=\!a_nS_n$ together with \eqref{LIL-MG-DR} immediately lead to
\begin{eqnarray*}
 \limsup_{n \rightarrow \infty} \Bigl(\frac{1}{2 n \log \log n}\Bigr)^{1/2} S_n & = & 
 -\liminf_{n \rightarrow \infty} \Bigl(\frac{1}{2 n \log \log n}\Bigr)^{1/2} S_n \nonumber \\
 & = & \frac{1}{\sqrt{3-4p}} \hspace{1cm} \text{a.s.}
\end{eqnarray*}
which completes the proof of Theorem \ref{T-ASCVGRATES-DR}.
\demend

\vspace{-5ex}
\subsection*{}
\begin{center}
{\bf B.2. The critical regime.}
\end{center}
\ \vspace{-4ex}\\

\noindent{\bf Proof of Theorem \ref{T-ASCVG-CR}.}
We shall proceed as in the proof of Theorem \ref{T-ASCVG-DR}. It follows from \eqref{CVGVNCR}
that, in the critical regime where $p=3/4$, $v_n$ increases slowly to infinity with a
logarithmic speed $\log n$. We obtain from Theorem 1.3.24 of \cite{Duflo97} that
\begin{equation}
\label{CVGRATEMNCR}
\frac{M_{n}^2}{\log n}  = O(\log \log n) \hspace{1cm} \text{a.s.}
\end{equation}
However, we deduce from \eqref{DEFAN} that
\begin{equation}
\label{CVGANCR}
\lim_{n \rightarrow \infty} n a_n^2 = \frac{\pi}{4}.
\end{equation}
Hence, as $M_n=a_nS_n$, we find from \eqref{CVGRATEMNCR} together with \eqref{CVGANCR} that
\begin{equation*}
\frac{S_{n}^2}{n\log n}  = O( \log \log n) \hspace{1cm} \text{a.s.}
\end{equation*}
which immediately implies \eqref{T-ASCVG-CR1}.
\demend
\ \vspace{-1ex}\\


\noindent{\bf Proof of Theorem \ref{T-ASCVGRATES-CR}.} The proof of Theorem
\ref{T-ASCVGRATES-CR} is left to the reader as it follows essentially the same lines as the one of
Theorem \ref{T-ASCVGRATES-DR}.
\demend
\ \vspace{-1ex}\\

\vspace{-5ex}
\subsection*{}
\begin{center}
{\bf B.3. The superdiffusive regime.}
\end{center}
\ \vspace{-4ex}\\

\noindent{\bf Proof of Theorem \ref{T-ASCVG-SR}.}
In the superdiffusive regime $3/4<p\leq 1$, we already saw from \eqref{CVGVNSR}
that $v_n$ converges to a finite value. Hence, as $\cMc_n \leq v_n$, we deduce 
from \eqref{DECMN} together with Theorem 1.3.15 of \cite{Duflo97}, the almost sure convergence
\begin{equation}
\label{CVGMNSR}
\lim_{n \rightarrow \infty} M_n  = M
\hspace{1cm} \text{where} \hspace{1cm} 
M=\sum_{k=1}^\infty a_k \veps_k.
\end{equation}
Moreover, we obtain from \eqref{DEFAN} that
\begin{equation}
\label{CVGANSR}
\lim_{n \rightarrow \infty} n^{2p-1} a_n = \Gamma(2p).
\end{equation}
Consequently, as $M_n=a_nS_n$, \eqref{T-ASCVG-SR1} immediately follows from
\eqref{CVGMNSR} and \eqref{CVGANSR}. One can observe that
$$
L=\frac{1}{\Gamma(2p)} \sum_{k=1}^\infty a_k \veps_k.
$$
It only remains to prove convergence \eqref{T-ASCVG-SR2}. For that purpose, it is only necessary to show that
the martingale $(M_n)$ is bounded in $\dL^4$.
It is not hard to see that $(M_n)$ is bounded in $\dL^2$. As a matter of fact, since $M_{n+1}^2=(a_{n+1} \veps_{n+1}+M_{n})^2$, 
it follows from \eqref{UBCEEPS} that
$$\dE[M_{n+1}^2| \cF_{n}]=a_{n+1}^2 \dE[\veps_{n+1}^2| \cF_{n}]+M_{n}^2 \leq
a_{n+1}^2 +M_{n}^2 \hspace{1cm} \text{a.s.}
$$
Taking expectation on both sides, we get that $\dE[M_{n+1}^2] \leq a_{n+1}^2 +\dE[M_{n}^2]$ leading to
\begin{equation}
\label{BL2}
\dE[M_n^2] \leq \sum_{k=1}^n a_k^2.
\end{equation}
Consequently, we obtain from \eqref{CVGVNSR} together with \eqref{BL2} that
\begin{equation}
\sup_{n \geq 1} \dE[M_n^2] \leq  {}_{3}F_2\Bigl( \begin{matrix}
{1,1,1}\\
{2p,2p}\end{matrix} \Bigl| 1\Bigr) < \infty.
\label{MBL2}
\end{equation}
By the same token, as $M_{n+1}^4=(a_{n+1} \veps_{n+1}+M_{n})^4$, we clearly have 
\begin{equation}
\dE[M_{n+1}^4 | \cF_{n}] = \sum_{\ell=0}^4 \binom{4}{\ell} a_{n+1}^\ell M_{n}^{4-\ell} \dE[\veps_{n+1}^\ell | \cF_{n}].
\label{CEM4}
\end{equation}
On the one hand, we already saw from \eqref{UBCEEPS} that 
\begin{equation*}
\dE[\veps_{n+1}^2 | \cF_{n}] \leq 1 \hspace{1cm} \text{and} \hspace{1cm}
\dE[\veps_{n+1}^4 | \cF_{n}] \leq \frac{4}{3} \hspace{1cm} \text{a.s.}
\end{equation*}
On the other hand, as in \eqref{CEEPS2}, we also have
\begin{eqnarray}
   \dE[\veps_{n+1}^3 | \cF_{n}]  
   & = &  2(\gamma_n -1)^3S_n^3  -2 (\gamma_n -1)^2 S_n  \hspace{1cm} \text{a.s.} \nonumber\\
   & = & 2(2p-1) \frac{S_n}{n} \Bigl( (2p-1)^2 \Bigl(\frac{S_n}{n}\Bigr)^2 -1 \Bigr)  \hspace{1cm} \text{a.s.}
   \label{CEEPS3}
\end{eqnarray}
Therefore, as $M_n=a_nS_n$ and $|S_n| \leq n$, we obtain from \eqref{CEEPS3} that
$$
M_n \dE[\veps_{n+1}^3 | \cF_{n}] \leq 0  \hspace{1cm} \text{a.s.}
$$
Consequently, it follows from \eqref{CEM4} that
$$
\dE[M_{n+1}^4 | \cF_{n}] \leq \frac{4}{3} a_{n+1}^4 + 6 a_{n+1}^2 M_n^2 + M_n^4
\hspace{1cm} \text{a.s.}
$$
Taking expectation on both sides, we find that
$$
\dE[M_{n+1}^4] \leq \frac{4}{3} a_{n+1}^4 + 6 a_{n+1}^2 \dE[M_n^2] + \dE[M_n^4]
$$
leading, via \eqref{BL2}, to
$$
\dE[M_n^4] \leq \frac{4}{3} \sum_{k=1}^n a_k^4  + 6  \sum_{k=1}^n a_k^2 \dE[M_{k-1}^2]
\leq 6\Bigl(1+\sum_{k=1}^n a_k^2 \Bigr) \sum_{k=1}^n a_k^2.
$$
Finally, we can deduce from \eqref{CVGVNSR} that
\begin{equation}
\sup_{n \geq 1} \dE[M_n^4] < \infty,
\label{MBL4}
\end{equation}
which completes the proof of Theorem \ref{T-ASCVG-SR}.
\demend
\ \vspace{-1ex}\\
The proof of Theorem \ref{T-MOM-SR} relies on the following well-known lemma on sums
of ratio of gamma functions.

\begin{lem}
\label{L-RG}
For any non-negative real numbers $a$ and $b$ such that $b\neq a+1$
and for all $n\geq 1$, we have
\begin{equation}
\label{SUMGAMMA}
\sum_{k=1}^n \frac{\Gamma(k+a)}{\Gamma(k+b)}= \frac{\Gamma(n+a+1)}{(b-a-1)\Gamma(n+b)}
\left(\frac{\Gamma(n+b) \Gamma(a+1)}{\Gamma(n+a+1)\Gamma(b)} -1\right).
\end{equation}
\end{lem}

\noindent{\bf Proof of Theorem \ref{T-MOM-SR}.}
Denote $\alpha=2p-1$, $\beta=2q-1$, 
$$
L_n=\frac{M_n}{\Gamma(2p)}
\hspace{1.5cm} \text{and} \hspace{1.5cm}
L=\frac{M}{\Gamma(2p)}.
$$
It follows from convergence \eqref{T-ASCVG-SR2} that for any integer $d=1,\ldots,4$
\begin{equation}
\label{LIMLD}
\lim_{n\rightarrow \infty} \dE[L_n^d]=\dE[L^d].
\end{equation}
First of all, we already saw that for all $n\geq 1$, $\dE[S_{n+1}|\cF_n]=\gamma_n S_n$ a.s. Consequently,
$$\dE[S_{n+1}]=\gamma_n \dE[S_n]=\Bigl(\frac{n+\alpha}{n}\Bigr)\dE[S_n],$$
which leads to
\begin{equation}
\label{MOMSN1}
\dE[S_n]=\prod_{k=1}^{n-1} \Bigl(\frac{k+\alpha}{k}\Bigr) \dE[S_1]= 
\frac{\beta \Gamma(n+\alpha)}{\Gamma(n) \Gamma(\alpha +1)}=
\frac{\beta}{a_n}.
\end{equation}
Hence, we immediately get from \eqref{MOMSN1} that
$$
\dE[L_n]=\dE[L]=
\frac{\beta}{\Gamma(\alpha +1)}
=\frac{2q-1}{\Gamma(2p)}.
$$
Next, taking expectation on both sides of \eqref{CES2}, we obtain that for all $n\geq 1$, 
$$
\dE[S_{n+1}^2] =1+ ( 2 \gamma_n -1)\dE[S_{n}^2]=1+ \Bigl(\frac{n+2\alpha}{n}\Bigr)\dE[S_n^2]
$$
which implies that 
\begin{eqnarray*}
\dE[S_n^2] & = &\frac{\Gamma(n+2\alpha)}{\Gamma(n)\Gamma(2\alpha +1)}\left(1+ \sum_{k=1}^{n-1} 
\frac{\Gamma(k+1)\Gamma(2\alpha +1)}{\Gamma(k+2\alpha+1)}\right), \\
& = &
\frac{\Gamma(n+2\alpha)}{\Gamma(n)}\sum_{k=1}^n \frac{\Gamma(k)}{\Gamma(k+2\alpha)}.
\end{eqnarray*}
Therefore, we deduce from identity \eqref{SUMGAMMA} with $a=0$ and $b=2\alpha$ that
\begin{equation}
\label{MOMSN2}
\dE[S_n^2]=\frac{n}{2\alpha -1}\left( \frac{\Gamma(n+2\alpha)}{\Gamma(n+1)\Gamma(2\alpha)} -1\right).
\end{equation}
Hence, we obtain from
\eqref{MOMSN2} that
\begin{equation*}
\dE[L_n^2]=\frac{a_n^2 \dE[S_n^2]}{(\Gamma(2p))^2}= \frac{n}{2\alpha -1} \left( \frac{\Gamma(n)}{\Gamma(n+\alpha)}\right)^2 \left( \frac{\Gamma(n+2\alpha)}{\Gamma(n+1)\Gamma(2\alpha)} -1\right)
\end{equation*}
which ensures, via \eqref{LIMLD} with $d=2$, that
$$
\lim_{n\rightarrow \infty} \dE[L_n^2]=\dE[L^2]=
\frac{1}{(2\alpha -1)\Gamma(2\alpha)}
=\frac{1}{(4p-3)\Gamma(2(2p-1))}.
$$
Furthermore, as $S_{n+1}^3=(S_n+X_{n+1})^3$, we obtain that
$$
\dE[S_{n+1}^3| \cF_n] = (3 \gamma_n -2)S_{n}^3 + (\gamma_n+2)S_n \hspace{1cm}\text{a.s}
$$
Consequently,
$$
\dE[S_{n+1}^3] =
\Bigl(\frac{3n+\alpha}{n}\Bigr)\dE[S_n]
+\Bigl(\frac{n+3\alpha}{n}\Bigr)\dE[S_n^3].
$$
Hence, it follows from tedious but straighforward calculations that
\begin{equation}
\label{MOMSN3SIGN}
\dE[S_n^3] =\frac{\beta\Gamma(n+3\alpha)}{\Gamma(n)\Gamma(3\alpha + 1)}\left(1+
\frac{3\Gamma(3\alpha +1)}{\Gamma(\alpha +1)}\xi_n\right)
\end{equation}
where
$$
\xi_n= \sum_{k=1}^{n-1} \Bigl(k+\frac{\alpha}{3}\Bigr)
\frac{\Gamma(k+\alpha)}{\Gamma(k+3\alpha+1)}
=\sum_{k=1}^{n-1} 
\frac{\Gamma(k+\alpha+1)}{\Gamma(k+3\alpha+1)}
-\frac{2\alpha}{3}\sum_{k=1}^{n-1} 
\frac{\Gamma(k+\alpha)}{\Gamma(k+3\alpha+1)}.
$$
However, we infer from \eqref{SUMGAMMA} with $a=\alpha$ or $a=\alpha+1$ and $b=3\alpha+1$ that
$$
\xi_n=\frac{1}{2\alpha-1} \left(
\frac{\Gamma(\alpha+2)}{\Gamma(3\alpha+1)} -
\frac{\Gamma(n+\alpha+1)}{\Gamma(n+3\alpha)}\right)
-\frac{1}{3}\left(
\frac{\Gamma(\alpha+1)}{\Gamma(3\alpha+1)} -
\frac{\Gamma(n+\alpha)}{\Gamma(n+3\alpha)}\right).
$$
Therefore, we obtain from \eqref{MOMSN3SIGN} that
\begin{eqnarray}
\dE[S_n^3] & = &\frac{\beta\Gamma(n+3\alpha)}{\Gamma(n)\Gamma(3\alpha + 1)}\left(
\frac{3(\alpha +1)}{2\alpha -1}
-\frac{\Gamma(3\alpha+1)\Gamma(n+\alpha)}{\Gamma(\alpha+1)\Gamma(n+3\alpha)}\frac{3n+\alpha +1}{2\alpha -1}
\right), \nonumber\\
& = & \frac{\beta}{(2\alpha -1)\Gamma(n)}\left(
\frac{3(\alpha +1)\Gamma(n+3\alpha)}{\Gamma(3\alpha + 1)}
-\frac{\Gamma(n+\alpha)}{\Gamma(\alpha+1)}(3n+\alpha +1)
\right).
\label{MOMSN3}
\end{eqnarray}
Consequently, we obtain from
\eqref{MOMSN3} that
\begin{equation*}
\dE[L_n^3]=
\frac{\beta(\Gamma(n))^2}{(2\alpha -1)(\Gamma(n+\alpha))^3} 
\left(
\frac{3(\alpha +1)\Gamma(n+3\alpha)}{\Gamma(3\alpha + 1)}
-\frac{\Gamma(n+\alpha)}{\Gamma(\alpha+1)}(3n+\alpha +1)
\right)
\end{equation*}
which leads, via \eqref{LIMLD} with $d=3$, to
$$
\lim_{n\rightarrow \infty} \dE[L_n^3]=\dE[L^3]=
\frac{\beta(\alpha+1)}{\alpha(2\alpha -1)\Gamma(3\alpha)}
=\frac{2p(2q-1)}{(2p-1)(4p-3)\Gamma(3(2p-1))}.
$$
By the same token, since $S_{n+1}^4=(S_n+X_{n+1})^4$, we obtain that
$$
\dE[S_{n+1}^4] =
1+2\Bigl(\frac{3n+2\alpha}{n}\Bigr)\dE[S_n^2]
+\Bigl(\frac{n+4\alpha}{n}\Bigr)\dE[S_n^4],
$$
which implies that
\begin{equation}
\label{MOMSN4PN}
\dE[S_n^4] =\frac{\Gamma(n+4\alpha)}{\Gamma(n)\Gamma(4\alpha + 1)}\left(1+
\frac{\Gamma(4\alpha +1)}{(2\alpha -1)\Gamma(2\alpha)}\bigl(P_n+6Q_n-6R_n\bigr)\right)
\end{equation}
where
\begin{eqnarray*}
P_n & = & (2\alpha -1)\Gamma(2\alpha) \sum_{k=1}^{n-1} 
\frac{\Gamma(k+1)}{\Gamma(k+4\alpha+1)}, \\
Q_n & = & \sum_{k=1}^{n-1}
\Bigl(k+\frac{2\alpha}{3}\Bigr) 
\frac{\Gamma(k+2\alpha)}{\Gamma(k+4\alpha+1)}, \\
R_n & = & \Gamma(2\alpha) \sum_{k=1}^{n-1}
\Bigl(k+\frac{2\alpha}{3}\Bigr) 
\frac{\Gamma(k+1)}{\Gamma(k+4\alpha+1)}.
\end{eqnarray*}
We make use once again of identity \eqref{SUMGAMMA} with appropriate values of
$a$ and $b$, to find that
\begin{equation}
\dE[S_n^4]   =  \frac{\Gamma(n+4\alpha)}{\Gamma(n)\Gamma(4\alpha + 1)}\left(
\frac{24\alpha(2\alpha(\alpha +1)-1)}{(2\alpha -1)^2 (4\alpha -1)}
 -  \frac{\Gamma(4\alpha+1)}{(2\alpha -1)^2 \Gamma(n+4\alpha)}\zeta_n \right)
\label{MOMSN4}
\end{equation}
where
$$
\zeta_n=\frac{2(3n+2(\alpha+1))}{\Gamma(2\alpha)}\Gamma(n+2\alpha)-
\frac{3n(4\alpha-1)+2(2\alpha^2+1)}{(4\alpha-1)}\Gamma(n+1).
$$
Finally, we deduce from \eqref{LIMLD} with $d=4$ that
\begin{eqnarray*}
\lim_{n\rightarrow \infty} \dE[L_n^4] & = &\dE[L^4]=
\frac{24\alpha(2\alpha(\alpha +1)-1)}{(2\alpha -1)^2 (4\alpha -1)\Gamma(4\alpha+1)}
 =  \frac{6(2\alpha(\alpha +1)-1)}{(2\alpha -1)^2 (4\alpha -1)\Gamma(4\alpha)}, \\
 & = &\frac{6(8p^2-4p-1)}{(8p-5)(4p-3)^2\Gamma(4(2p-1))},
\end{eqnarray*}
 which completes the proof of Theorem \ref{T-MOM-SR}.
\demend
\vspace{-2ex}

\section*{Appendix C \\ Proofs of the asymptotic normality results}
\renewcommand{\thesection}{\Alph{section}}
\renewcommand{\theequation}{\thesection.\arabic{equation}}
\setcounter{section}{3}
\setcounter{equation}{0}

\vspace{-2ex}
\subsection*{}
\begin{center}
{\bf C.1. The diffusive regime.}
\end{center}
\ \vspace{-4ex}\\

\noindent{\bf Proof of Theorem \ref{T-AN-DR}.}
We shall make use of the central limit theorem for martingales given e.g. by Corollary 2.1.10 of  \cite{Duflo97},
to establish the asymptotic normality \eqref{T-AN-DR1}. It follows from \eqref{T-ASCVG-DR1},
\eqref{CALIP} and \eqref{CVGVNDR} that
\begin{equation}
\label{CVG-IP-DR}
\lim_{n\rightarrow \infty} 
\frac{\cMc_{n}}{v_n} = 1  \hspace{1cm} \text{a.s.}
\end{equation}
Hereafter, it only remains to prove that $(M_n)$ satisfies Lindeberg's condition, that is, for all $\veps > 0$, 
\begin{equation}
\label{LINDEBERG}
\frac{1}{v_n} \sum_{k=1}^{n} \dE \Bigl[ | \Delta M_{k} |^2 \rI_{| \Delta M_{k}|    
\geq \veps \sqrt{v_n}} \bigl| \cF_{k-1} \Bigr] \limp 0.
\end{equation}
We obtain from \eqref{UBCEEPS} that for all $\veps > 0$,
\begin{eqnarray*}
\frac{1}{v_n} \sum_{k=1}^{n} \dE \Bigl[ | \Delta M_{k} |^2 \rI_{| \Delta M_{k}| \geq \veps \sqrt{v_n}} \bigl| \cF_{k-1} \Bigr] 
& \leq & \frac{1}{\veps^2v_n^2} \sum_{k=1}^{n} \dE \left[ | \Delta M_{k} |^4  | \cF_{k-1} \right],   \\
& \leq & \sup_{1 \leq k \leq n } \dE[\veps_{k}^4 | \cF_{k-1}] \frac{1}{\veps^2v_n^2} \sum_{k=1}^{n} a_k^4,   \\
& \leq & \frac{4}{3\,\veps^2v_n^2} \sum_{k=1}^{n} a_k^4.
\end{eqnarray*}
However, we deduce from \eqref{CONDLIL} together with Kronecker's lemma that
\begin{equation*}
\lim_{n \rightarrow \infty} \frac{1}{v_n^2} \sum_{k=1}^{n} a_{k}^4=0,
\end{equation*}
which ensures that Lindeberg's condition is satisfied. Therefore, we can conclude from the central limit theorem 
for martingales that
\begin{equation}
\label{CLTMN-DR}
\frac{1}{\sqrt{v_n}} M_{n} \liml \cN (0, 1).
\end{equation}
As $M_n=a_n S_n$ and $\sqrt{n}a_n$ is equivalent to $\sqrt{v_n(3-4p)}$, we find from
\eqref{CLTMN-DR} that
\begin{equation*}
\frac{S_n}{\sqrt{n}} \liml \cN \Bigl(0, \frac{1}{3-4p} \Bigr),
\end{equation*}
which is exactly what we wanted to prove.
\demend

\vspace{-5ex}
\subsection*{}
\begin{center}
{\bf C.2. The critical regime.}
\end{center}
\ \vspace{-4ex}\\

\noindent{\bf Proof of Theorem \ref{T-AN-CR}.}
We shall now proceed to the proof of the asymptotic normality \eqref{T-AN-CR1}
making use, once again, of the central limit theorem for martingales given e.g. by Corollary 2.1.10 of  \cite{Duflo97}.
We have from \eqref{T-ASCVG-CR1}, \eqref{CALIP} and \eqref{CVGVNCR} that
\begin{equation}
\label{CVG-IP-CR}
\lim_{n\rightarrow \infty} 
\frac{\cMc_{n}}{v_n} = 1  \hspace{1cm} \text{a.s.}
\end{equation}
Moreover, it follows from \eqref{CVGVNCR} and \eqref{CVGANCR} that
$a_n^2v_n^{-1}$ is equivalent to $(n\log n)^{-1}$. However, it is well-known that
\begin{equation*}
\sum_{n=1}^{+\infty} \frac{1}{n^2(\log n)^2} < \infty.
\end{equation*}
Consequently, we clearly have
\begin{equation}
\label{CONDLILCR}
\sum_{n=1}^{+\infty} \frac{a_{n}^4}{v_n^2} < +\infty. 
\end{equation}
We already saw in the proof of Theorem \ref{T-AN-DR} that \eqref{CONDLILCR} is enough
to check that $(M_n)$ satisfies Lindeberg's condition. Consequently, we can deduce from the 
central limit theorem for martingales that
\begin{equation}
\label{CLTMN-CR}
\frac{1}{\sqrt{v_n}} M_{n} \liml \cN (0, 1).
\end{equation}
Hence, as $M_n=a_n S_n$ and $\sqrt{n \log n}a_n$ is equivalent to $\sqrt{v_n}$, we find from
\eqref{CLTMN-CR} that
\begin{equation*}
\frac{S_n}{\sqrt{n \log n}} \liml \cN (0,1),
\end{equation*}
which achieves the proof of Theorem \ref{T-AN-CR}.
\demend

\vspace{-4ex}

\bibliographystyle{acm}

\begin{thebibliography}{1}

\bibitem{Baur16}
{\sc Baur, E. and Bertoin, J.} 
\newblock Elephant Random Walks and their connection to P\'olya-type urns.
\newblock {\em Phys. Rev. E 94}, 052134 (2016).

\bibitem{Boyer14}
{\sc Boyer, D., Romo-Cruz, J. C. R.} 
\newblock Solvable random-walk model with memory and its relations
with Markovian models of anomalous diffusion. 
\newblock {\em Phys. Rev. E 90}, 042136 (2014).

\bibitem{Bercu04}
{\sc Bercu, B.}
\newblock On the convergence of moments in the almost sure central limit theorem for martingales 
with statistical applications. 
\newblock {\em Stochastic Process. Appl. 111}, 1 (2004), pp. 157--173.

\bibitem{Chauvin11}
{\sc Chauvin, B., Pouyanne, N., Sahnoun, R.} 
\newblock Limit distributions for large P\'olya urns. 
\newblock {\em Ann. Appl. Probab. 21}, (2011), pp 1–-32.

\bibitem{Col17}
{\sc Coletti, C. F., Gava, R., Sch\"utz, G. M.} 
\newblock Central limit theorem and related results for the elephant random walk.
\newblock {\em J. Math. Phys. 58}, 053303 (2017).

\bibitem{Cressoni13}
{\sc Cressoni, J.~C., Viswanathan, G.~M., Da Silva, M.~A.~A., }
\newblock Exact solution of an anisotropic 2D random walk
model with strong memory correlations.
\newblock {\em J. Phys. A: Math. Theor. 46}, 505002 (2013).

\bibitem{Cressoni07}
{\sc Cressoni, J.~C., Da Silva, M.~A.~A., Viswanathan, G.~M.}
\newblock Amnestically induced persistence in random walks. 
\newblock {\em Phys. Rev. Let. 98}, 070603 (2007).


\bibitem{Da13}
{\sc Da Silva, M.~A.~A., Cressoni, J.~C., Sch\"utz, G.~M., Viswanathan, G.~M., Trimper, S.}
\newblock Non-Gaussian propagator for elephant random walks. 
\newblock {\em Phys. Rev. E 88}, 022115 (2013).

\bibitem{Duflo97}
{\sc Duflo, M.}, 
\newblock Random iterative models, Vol. 34 of Applications of Mathematics. 
Springer-Verlag, Berlin, 1997.

\bibitem{Janson04}
\newblock Functional limit theorems for multitype branching processes and generalized P\'{o}lya urns. 
\newblock {\em Stochastic Process. Appl. 110}, 2 (2004), pp. 177--245.

\bibitem{Kumar10}
{\sc Kumar, N., Harbola, U., Lindenberg, K.}
\newblock Memory-induced anomalous dynamics: Emergence
of diffusion, subdiffusion, and superdiffusion from a single random walk model.
\newblock {\em Phys. Rev. E 82}, 021101 (2010).

\bibitem{Harris15}
{\sc Harris, R.}
\newblock Random walkers with extreme value memory: modelling the peak-end rule.
\newblock {\em New J. Phys. 17}, 053049 (2015).

\bibitem{HallHeyde80}
{\sc Hall, P., and Heyde, C. C.} 
\newblock Martingale limit theory and its application.
Academic Press Inc., New York,  1980.

\bibitem{Harbola14}
{\sc Harbola, U., Kumar, N., Lindenberg, K.}
\newblock Memory-induced anomalous dynamics in a minimal random walk model.
\newblock {\em Phys. Rev. E 90}, 022136 (2014).


\bibitem{Kursten16}
{\sc K\"ursten, R.}
\newblock Random recursive trees and the elephant random walk. 
\newblock {\em Phys. Rev. E 93}, 032111 (2016).

\bibitem{Lyu17}
{\sc Lyu, J., Xin, J., Yu, Y.}
\newblock Residual diffusivity in elephant random
walk models with stops. 
\newblock {\em arXiv:1705.0271}, (2017).


\bibitem{Par06}
{\sc Paraan, F. N. C. and Esguerra, J. P.}
\newblock Exact moments in a continuous time random walk with complete memory of its history. 
\newblock {\em Phys. Rev. E 74}, 032101 (2006).

\bibitem{Schutz04}
{\sc Sch\"utz, G.~M., and Trimper, S.}
\newblock Elephants can always remember: Exact long-range memory
effects in a non-Markovian random walk. 
\newblock {\em Phys. Rev. E 70}, 045101 (2004).

\bibitem{Stout70}
{\sc Stout, W. F.}
\newblock A martingale analogue of Kolmogorov’s law of the iterated logarithm.
\newblock {\em Z. Wahrscheinlichkeitstheorie 15} (1970), pp. 279–-290.

\bibitem{Stout74}
{\sc Stout, W.~F.}, 
\newblock Almost sure convergence, Probability and Mathematical Statistics, Vol. 24, 
Academic Press, New York-London, 1974.



\end{thebibliography}

\end{document}